\documentclass[10pt]{amsart}
\usepackage{amssymb}
\usepackage{amsfonts}
\usepackage{amssymb,latexsym}
\usepackage{enumerate}
\usepackage{mathrsfs}
\usepackage{amssymb}
\makeatletter
\@namedef{subjclassname@2010}{%
  \textup{2010} Mathematics Subject Classification}
\makeatother

  [2002/01/19 v2.2g %
    AMS font definitions%
  ]
\DeclareFontFamily{U}{euf}{}
\DeclareFontShape{U}{euf}{m}{n}{%
  <5><6><7><8><9>gen*eufm%
  <10><10.95><12><14.4><17.28><20.74><24.88>eufm10%
  }{}
\DeclareFontShape{U}{euf}{b}{n}{%
  <5><6><7><8><9>gen*eufb%
  <10><10.95><12><14.4><17.28><20.74><24.88>eufb10%
  }{}

  [2002/01/19 v2.2g %
    AMS font definitions%
  ]
\DeclareFontFamily{U}{msb}{}
\DeclareFontShape{U}{msb}{m}{n}{%
  <5><6><7><8><9>gen*msbm%
  <10><10.95><12><14.4><17.28><20.74><24.88>msbm10%
  }{}

  [2002/01/19 v2.2g %
    AMS font definitions%
  ]
\DeclareFontFamily{U}{msa}{}
\DeclareFontShape{U}{msa}{m}{n}{%
  <5><6><7><8><9>gen*msam%
  <10><10.95><12><14.4><17.28><20.74><24.88>msam10%
  }{}

\newtheorem{theorem}{Theorem}[section]
\newtheorem{lemma}[theorem]{Lemma}

\theoremstyle{definition}
\newtheorem{remark}[theorem]{Remark}

\newcommand{\bell}{\textup{B}}

\numberwithin{equation}{section} 

\def\C{\mathbb C_p}
\def\BZ{\mathbb Z}
\def\Z{\mathbb Z_p}
\def\Q{\mathbb Q_p}
\def\C{\mathbb C_p}
\def\BZ{\mathbb Z}
\def\Z{\mathbb Z_p}
\def\Q{\mathbb Q_p}

\DeclareMathOperator{\td}{d}
\begin{document}

\title[]
{Two closed forms for the Apostol-Bernoulli polynomials}

\author{Su Hu}
\address{Department of Mathematics, South China University of Technology, Guangzhou, Guangdong 510640, China}
\email{mahusu@scut.edu.cn}

\author{Min-Soo Kim}
\address{Center for General Education, Kyungnam University,
7(Woryeong-dong) kyungnamdaehak-ro, Masanhappo-gu, Changwon-si, Gyeongsangnam-do 631-701, Republic of Korea}
\email{mskim@kyungnam.ac.kr}


\begin{abstract}
In this note, we shall obtain two closed forms for the Apostol-Bernoulli polynomials.
\end{abstract}



\maketitle

\def\B{\overline B}
\def\E{\overline E}
\def\Z{\mathbb Z_p}



\def\ord{\text{ord}_p}
\def\C{\mathbb C_p}
\def\BZ{\mathbb Z}
\def\Z{\mathbb Z_p}
\def\Q{\mathbb Q_p}
\def\wh{\widehat}
\def\ov{\overline}


\section{Introduction}
\label{Intro}
The Bernoulli numbers $B_{n}$ are defined by the generating function
\begin{equation}\label{Bn}
\frac{x}{e^x-1}=\sum_{n=0}^\infty B_n\frac{x^n}{n!}=1-\frac x2 +\sum_{n=1}^\infty B_{2n}\frac{x^{2n}}{(2n)!},\quad |x|<2 \pi
\end{equation}
and the Bernoulli polynomials $B_{n}(u)$ are defined by
\begin{equation}\label{Bp}
\frac{xe^{ux}}{e^x-1}=\sum_{n=0}^{\infty}B_{n}(u) \frac{x^n}{n!},\quad |x|<2 \pi.
\end{equation}
Because the function $\frac{x}{e^x-1}-1+\frac x2$ is an even function, all of the Bernoulli numbers
with odd subscripts $>1$ are zero.
It is clear that $B_0=1$ and $B_1=-\frac 12$ (see \cite{Qi}).

The Apostol-Bernoulli polynomials $B_n(u,z)$ are natural generalizations of the Bernoulli polynomials, they were  introduced  by
Apostol in 1951 (see \cite{Ap}) and were first named by Luo  in 2004 (see \cite{Luo04}). Their definitions are as follows,
\begin{equation}\label{AB-def}
\left(\frac{x}{z e^x-1}\right) e^{ux}=\sum_{n=0}^\infty B_n(u,z) \frac{x^n}{n!},
\end{equation}
where $|x|< 2\pi$ when $z=1$; $|x|< |\textrm{log}~z|$ when $z\neq 1$ (see ~\cite{Luo04}). In particular, $B_n(z)=B_n(0,z)$ are the
Apostol-Bernoulli numbers. Letting $z=1$ in (\ref{AB-def}), we
obtain the Bernoulli polynomials $B_n(u)$ and Bernoulli
numbers $B_n,$ respectively.

As pointed out by Carlitz (see \cite{Ap2}), formula (3.3) of \cite{Ap} leads to the result
\begin{equation*}
z^{m}B_{n}(m,z)-B_{n}(0,z)=n\sum_{a=0}^{m-1}a^{n-1}z^{a}
\end{equation*}
which, for integer values of the variable $u$, gives a relation of the functions $B_{n}(u,z)$ with the Mirimanoff polynomials discussed by Vandiver in \cite[p.504]{Van}.
The Mirimanoff polynomials are written in the form
\begin{equation}\label{Miri}
f_{n}^{(m)}(z)=0^{n}+1^{n}z+2^{n}z^{2}+\cdots+(m-1)^{n}z^{m-1},
 \end{equation}
 where $0^{0}=1$ (see \cite[p.~504, Eq.~(1)]{Van}). In \cite{Van}, Vandiver also pointed out that the Apostol-Bernoulli polynomials may also be tracked back
 to L. Euler's classical work ``Institutiones calculi differentialis, (II) 487--491" published in 1755 \cite{Euler}.
 In fact, if we write
 \begin{equation}
 H_{n}=R_{n}(z)/(z-1)^{n},
 \end{equation}
 here $R_{n}(z)$ are the Euler polynomials defined by
 \begin{equation}
 \begin{aligned}
 R_{1}(z)&=1,\\
 R_{2}(z)&=1+z,\\
 R_{3}(z)&=1+4z+z^{2},\\
 R_{4}(z)&=1+11z+11z^{2}+z^{3},\\
 R_{5}(z)&=1+26z+66z^{2}+26z^{3}+z^{4},\\
  \cdots,
  \end{aligned}
 \end{equation}
 then we have \cite[p.~508, Eq.~(19)]{Van}
 \begin{equation}
z^{r+1}[(H-r)^{n}-zH_{n}=(-1)^{n}(z-1)f_{n}^{(r+1)}(z),
 \end{equation}
 where $f_{n}^{(r+1)}(z)$ are the Mirimanoff polynomials (see Eq.~(\ref{Miri}) above).
 For the notations and the details, we refer to Vandiver in \cite[p.~507-508]{Van}.

The  Apostol-Bernoulli polynomials and numbers including their applications have been widely studied
  by many authors.  In 1950, Apostol~\cite{Ap} presented may fundamental properties the Apostol-Bernoulli polynomials $B_n(u,z)$, and obtained the recursion formula for the Apostol-Bernoulli numbers $B_{n}(z)$ using the Stirling numbers of the
  second kind (see \cite[p. 166, Eq.~(3.7)]{Ap} as follows
  \begin{equation}\label{ApR}
  B_{n}(z)=n\sum_{k=1}^{n-1}(-1)^{k}k!z^{k}(z-1)^{-1-k}S(n-1,k),\quad (n\in\mathbb{N}_{0};\, \textrm{Re}(z)>0;\, z\neq 1)
  \end{equation} where $S(n,k)$ denote the Stirling numbers of the second kind which are defined by means of the following expressions
  (see \cite[p.~207, Theorem B]{Comtet-Combinatorics-74}) \begin{equation} x^{n}=\sum_{k=0}^{n}\binom{x}{k}k!S(n,k).\end{equation} In 2000, by applying Lerch's functional equation, Srivastava \cite[p.~84, Eq.~(4.6)]{Srivastava2000}  obtained the following representations of the Apostol-Bernoulli polynomials $B_{n}(u,z)$ at rational points $u=\frac{p}{q}$ in terms of the Hurwitz zeta function \begin{equation}\label{Srivastava}\begin{aligned} B_{n}\left(\frac{p}{q}, e^{2\pi i\xi}\right)&=-\frac{n!}{(2q\pi)^{n}}\Bigg\{\sum_{j=1}^{q}\zeta\bigg(n,\frac{\xi+j-1}{q}\bigg)\textrm{exp}\bigg[\bigg(\frac{n}{2}-\frac{2(\xi+j-1)p}{q}\bigg)\pi i\bigg]\\
 &\quad+\sum_{j=1}^{q}\zeta\bigg(n,\frac{j-\xi}{q}\bigg)\textrm{exp}\bigg[\bigg(-\frac{n}{2}+\frac{2(j-\xi)p}{q}\bigg)\pi i\bigg]\Bigg\},
 \end{aligned}\end{equation}
 where $n\in\mathbb{N}\setminus\{1\},q\in\mathbb{N}, p\in\mathbb{Z}$ and $\xi\in\mathbb{R}.$
 In 2004, by applying binomial series expansion and Leibniz's rule, Luo   obtained the following representation of the Apostol-Bernoulli  polynomials $B_{n}(u,z)$ by the Gaussian hypergeometric functions \cite[p. 510, Theorem 2.1]{Luo04}
 \begin{equation}\label{Luo2004Eq}\begin{aligned} B_{n}(u,z)&=n\sum_{l=0}^{n-1}\binom{n-1}{l}z^{l}(z-1)^{-l-1}\sum_{j=0}^{l}(-1)^{j}\binom{l}{j}j^{l}(u+j)^{n-l-1}\\
  &\quad\times  \sideset{{}_2}{_1}{\mathop{F}}(l-n-1,l; l+1; j/(u+j)),\end{aligned}\end{equation}
where \begin{equation} \sideset{{}_2}{_1}{\mathop{F}}(a,b;c;z):=\sum_{n=0}^{\infty}\frac{(a)_{n}(b)_{n}}{(c)_{n}}\frac{z^{n}}{n!}, \quad |z|<1\end{equation}
and $(z)_{n}$ denotes the Pochhammer symbol which is defined by
$(z)_{0}=1, (z)_{n}=z(z+1)\cdots(z+n-1)=\frac{\Gamma(z+n)}{\Gamma(z)}, (n\in\mathbb N).$ Luo's formula (\ref{Luo2004Eq}) not only generalized Apostol's formula (\ref{ApR}) (see \cite[p.~512, Remark 2.6]{Luo04}), but also led to the first known closed form for the Apostol-Bernoulli polynomials as follows \begin{equation}\label{LuoClosed} B_n(u,z)=\sum_{k=1}^n k \binom nk \sum_{j=0}^{k-1} (-1)^jz^j(z-1)^{-j-1}j! S(k-1,j)u^{n-k}, \quad(z\neq 1)\end{equation} (see \cite[p.~512, Remark 2.6]{Luo04}).

According to Wiki~\cite{Wiki},  ``In mathematics, a closed-form expression is a mathematical expression that can be evaluated in a finite number of operations. It may contain constants, variables, certain `well-known' operations (e.g., $+~-~\times ~\div$), and functions (e.g., $n$th root, exponent, logarithm, trigonometric functions, and inverse hyperbolic functions), but usually no limit."

Luo \cite{Luof}, Bayad \cite{Ba}, Navas, Francisco and Varona \cite{NRV} investigated Fourier expansions for the Apostol-Bernoulli and Apostol-Euler polynomials.
Grag, Jain and Srivastava \cite{GJS} investigated the generalized Apostol-Bernoulli polynomials of higher order, which were introduced and studied by Luo and Srivastava in \cite{Luo2006JMAA, LuoSrivastava2006}. They derived an explicit representation of these generalized Apostol-Bernoulli polynomials of higher order in terms of a generalization of the Hurwitz-Lerch zeta function, and  also established a functional relationship between the generalized Apostol-Bernoulli polynomials of rational arguments and the Hurwitz (or generalized) zeta function. Their results  extended the classical ones by  Apostol \cite{Ap} and Srivastava \cite{Srivastava2000}. In~\cite{CJS}, Choi, Jang and Srivastava presented an explicit representation of the generalized Bernoulli polynomials in terms of a generalization for the Hurwitz-Lerch zeta function. In \cite{Luo2006TWJ}, Luo defined the Apostol-Euler numbers and polynomials of higher order. He also established many of  their fundamental properties, including several explicit formulas involving the Gaussian hypergeometric function and the Stirling numbers of the second kind. He also deduced their special cases and applications that lead to the corresponding formulas of the classical Euler numbers and polynomials of higher order.  In \cite{LuoSrivastava2006}, Luo and Srivastava  derived several general properties and relationships involving the Apostol-Bernoulli and Apostol-Euler polynomials. Their work generalized the properties and relationships involving the classical as well as the generalized (or higher-order) Bernoulli and Euler polynomials obtained by Srivastava and Pint\'er in \cite{SrivastavaPinter} and Cheon in \cite{Cheon}. Furthermore, they also obtained several relationships associated with the Stirling numbers of the second kind.  In 2009, Luo \cite[p. 380, Theorem 2.1]{Luo2009IT} obtained the following multiplication formulas for the Apostol-Bernoulli polynomials of higher order
\begin{equation}\label{Luo2009Th2.1} B_{n}^{(\alpha)}(mx,z)=m^{n-\alpha}\sum_{v_{1},v_{2},\ldots, v_{m-1}\geq 0}\binom{\alpha}{v_{1},v_{2},\ldots, v_{m-1}}z^{r}B_{n}^{(\alpha)}\left(x+\frac{r}{m},z^{m}\right),\end{equation} where $r=v_{1}+2v_{2}+\cdots+(m-1)v_{m-1}.$ This generalized Carlitz's  multiplication formulas for the Bernoulli polynomials of higher order (see \cite[p.~381, Corollary 2.3]{Luo2009IT}). Let \begin{equation}S_{n}^{(\ell)}(m;\lambda)=\sum_{\substack{0\leq v_{1},v_{2},\ldots,v_{m}\leq
\ell\\v_{1}+v_{2}+\cdots+v_{m}=\ell}}\binom{\ell}{v_{1},v_{2},\cdots,v_{m}}\lambda^{v_{1}+2v_{2}+\cdots+mv_{m}}(v_{1}+2v_{2}+\cdots+mv_{m})^{n}\end{equation} be the $\lambda$-multiple power sums. In the same paper, Luo also gave the following representations of $S_{n}^{(\ell)}(m;\lambda)$ in terms of  the Apostol-Bernoulli polynomials and numbers of higher order  \cite[p. 382, Theorem 2.2]{Luo2009IT}  \begin{equation}\label{Luo2009Th2} S_{n}^{(\ell)}(m;\lambda)=\sum_{j=0}^{\ell}\binom{\ell}{j}\frac{(-1)^{\ell-j}\lambda^{mj+\ell}}{(n+1)_{\ell}}\sum_{k=0}^{n+\ell}\binom{n+\ell}{k}B_{k}^{(j)}(mj+\ell,\lambda)B_{n+\ell-k}^{(\ell-j)}(\lambda),\end{equation}  where $(n)_{0}=1, (n)_{k}=n(n+1)\cdots(n+k-1)$. This generalized the  following formula on the representation of power sums in terms of the Bernoulli polynomials which obtained by J. Bernoulli in 1713 \begin{equation} \sum_{j=1}^{m}j^{n}=\frac{B_{n+1}(m+1)-B_{n+1}}{n+1}\end{equation} (see \cite[p.~383, Corollary 2.6]{Luo2009IT}). Furthermore, Luo also proved the following recursive formula for the Apostol-Bernoulli numbers of higher order in terms of the $\lambda$-multiple power sums \cite[p. 384, Theorem 2.3]{Luo2009IT}\begin{equation}\label{Luo2009Th2.3}
m^{\ell}B_{n}^{(\ell)}(\lambda)-m^{n}B_{n}^{(\ell)}(\lambda^{m})=\sum_{k=0}^{n}\binom{n}{k}m^{k}B_{k}^{(\ell)}(\lambda^{m})S_{n-k}^{(\ell)}(m-1;\lambda).\end{equation}
This generalized the following recursive formula for the  Bernoulli numbers in terms of the   power sums which was obtained by  Howard  \cite[p.~160, Eq.(6)]{Howard}  or Deeba and Rodriguez \cite[p.~424, Eq.(8)]{Deeba}
\begin{equation}
(m-m^{n})B_{n}=\sum_{k=0}^{n}\binom{n}{k}m^{k}B_{k}S_{n-k}(m-1)
\end{equation} (see \cite[p.~384, Corollaries 2.7 and 2.8]{Luo2009IT} and the remarks below).
In \cite{Luo2009IT}, the parallel  results on the Apostol-Euler numbers and polynomials of higher order and the $\lambda$-multiple alternating sums have also been presented (see \cite[Theorems 3,1, 3.2 and 3.3]{Luo2009IT}). For example,  \cite[Theorem 3.2]{Luo2009IT} extended the following Euler's formula on the representation of alternating power sums in terms of the Euler polynomials  \begin{equation} \sum_{j=1}^{m}(-1)^{j+1}j^{n}=-\frac{(-1)^{m}E_{n}(m+1)+E_{n}(0)}{2}\end{equation} (see \cite[p.~388, Corollary 3.7]{Luo2009IT}). In 2010, Luo~\cite{LuoHou2010} introduced the notions of $\lambda$-Stirling numbers of the second kind and gave many of its fundamental properties. He also obtained an explicit relationship between the generalized Apostol-Bernoulli and Apostol-Euler polynomials in terms of the $\lambda$-Stirling numbers of the second kind.  In 2011, Srivastava \cite{Srivastava2011} gave  a survey to many results on the generalizations of the classical Euler, Bernoulli, and  Genocchi polynomials.  By collecting and studying the generating functions for each class of polynomials, he also gave a unified exposition of relations between the polynomials and highlighted their connections with relevant special functions.
Using the  $\lambda$-Stirling numbers of the second kind, he derived several properties and formulas and considered some of their interesting applications to the family of the Apostol-type polynomials. He also gave  a brief expository and historical account of the various basic (or $q$-) extensions of the classical Bernoulli polynomials and numbers, the classical Euler polynomials and numbers, the classical Genocchi polynomials and numbers, and also of their generalizations to the above-mentioned families of the Apostol-type polynomials and numbers. Furthermore, he also indicated many relevant connections of the definitions and results presented in this survey with those in earlier as well as forthcoming investigations. In \cite{LuoSrivastava2011}, Luo and Srivastava introduced the notions of the Apostol-Genocchi polynomials of higher order. They also established many fundamental properties of these polynomials, including  some explicit relationships with the Apostol-Bernoulli  and  Apostol-Euler polynomials. They derived some explicit series representations of these polynomials in terms of the Gaussian hypergeometric function and the Hurwitz (or generalized) zeta functions. These results generalized the corresponding ones for the Genocchi and Euler polynomials of higher order. Furthermore, using  the $\lambda$-Stirling numbers of the second kind, they also obtained several  basic properties and formulas related to the family of the Apostol-type polynomials.  In \cite{LuoOsaka2011}, Luo derived the Fourier expansions and integral representations for the Apostol-Genocchi polynomials. He  obtained their formulas at rational arguments in terms of Hurwitz zeta function and showed an explicit relationship with Gaussian hypergeometric functions. These generalized many fundamental results on the classical Genocchi polynomials. In particular, the Euler's famous formula \begin{equation}\zeta(2n)=\frac{(-1)^{n-1}(2\pi)^{2n}}{2(2n)!}B_{2n}\end{equation} has been derived in a different way. In 2012, Kim and Hu \cite{KH} obtained the sums of products identity for the Apostol-Bernoulli numbers which is an analogue of the classical sums of products identity for Bernoulli numbers dating back to Euler. In \cite{Srivastava2},  Srivastava, Kurt and Simsek constructed the generating functions for several families of Genocchi type polynomials, they defined a function which interpolates these polynomials at negative integers by applying the derivative operator to these generating functions, they proved a multiplication theorem for these polynomials, they also proved several other identities and provided many applications associated with these and related polynomials and their interpolation functions. In 2013, Lu and Luo \cite{LuLuo2013} obtained several  new properties of the generalized Apostol-type polynomials, including the recurrence relations, the differential equations and some other connected problems, which extended some known results. In the same paper, many new properties on the generalized Apostol-Euler polynomials, the generalized Apostol-Bernoulli polynomials, and Apostol-Genocchi polynomials of high order have also been derived. In 2014, Luo~\cite{LuoFilo} obtained the generating functions and basic properties of the $q$-Apostol-Bernoulli and $q$-Apostol-Euler polynomials. He also derived some relationships between the $q$-Apostol-Bernoulli and $q$-Apostol-Euler polynomials which generalized   some known results on  the corresponding $q$-Bernoulli and $q$-Euler polynomials.  In the same paper, several formulas in series of $q$-Stirling numbers of the second kind have  also been derived. In \cite{LuoUkr},  Luo obtained some further properties for the  $q$-Apostol-Euler polynomials, including Raabe's multiplication theorem and alternating sums representations for these polynomials. He also derived the $q$-extensions of some formulas in \cite{Luo2009IT}. Furthermore, he also obtained several formulas and relationships between the $q$-Apostol-Euler polynomials, $q$-Apostol-Euler polynomials, and $q$-Hurwitz-Lerch zeta functions. In 2015, Luo~\cite{Luo2015AMC} introduced the elliptic analogues of the Apostol-Bernoulli and Apostol-Euler polynomials. He also obtained  the closed expressions of sums of products for these elliptic type polynomials. In \cite{LuLuo2015}, Lu and Luo obtained
further properties of  the Apostol-type polynomials, including a unified multiplication formula  of these polynomials, the explicit representations of these polynomials in terms of the Gaussian hypergeometric function and the generalized Hurwitz zeta function.  Their results extended several previous known results by Luo, Garg, Srivastava, Ozden, and \"Ozarslan etc.

Recently, Qi and Chapman~\cite{Qi} got two closed forms for Bernoulli polynomials. In this note, we show that, different with~\cite{Qi}, if directly applying the  generating functions instead of their integral expressions, then Qi and Chapman's results \cite[Theorems 1.1 and 1.2]{Qi} may be generalized to other special functions. As a result,  in addition to Luo's formula (\ref{LuoClosed}) above, we  shall get the following two other closed forms for Apostol-Bernoulli polynomials and numbers, which may be used to compute these special polynomials and numbers in a finite number of steps.

\begin{theorem}\label{th-main1}
Suppose that $z\not=1$. The Apostol-Bernoulli polynomials $B_n(u,z)$ for $n\in\mathbb N$ may be expressed as
\begin{equation}\label{AP-Stirl-form}
\begin{aligned}
B_n(u,z)
&=n\sum_{k=0}^{n-1}\frac{(-1)^{k}k!}{(z-1)^{k+1}}\sum_{r+s=k}\sum_{\ell+m=n-1}(-1)^{s+m}\binom{n-1}{\ell}\\
&\quad\times z^{r}(1-u)^{\ell}u^m S(\ell,r) S(m,s).
 \end{aligned}
\end{equation}
 Consequently, the Apostol-Bernoulli numbers $B_n(z)$ for $n\in\mathbb N$ can be represented as
\begin{equation}\label{APN-Stirl-form}
\begin{aligned}
B_n(z)=n\sum_{k=0}^{n-1}\frac{(-1)^{k}k!}{(z-1)^{k+1}}z^{k}S(n-1,k).
 \end{aligned}
\end{equation}
 \end{theorem}
 \begin{remark}Suppose that $z\not=1$. Letting $x=0$ in the both sides of (\ref{AB-def}), we get $B_{0}(u,z)=0$ and $B_{0}(z)=0.$
 Since by (\ref{AB-def}), we have
 $$ \frac{d}{dx}\left(\frac{xe^{ux}}{ze^{x}-1}\right)\biggl|_{x=0}=\sum_{m=1}^{\infty}B_{m}(u,z)\frac{x^{m-1}}{(m-1)!}\biggl|_{x=0}=B_{1}(u,z)$$
 and \begin{equation}\begin{aligned}
 \frac{d}{dx}\left(\frac{xe^{ux}}{ze^{x}-1}\right)&=\frac{e^{ux}}{ze^{x}-1}+x\frac{d}{dx}\left(\frac{e^{ux}}{ze^{x}-1}\right)\\
&\rightarrow \frac{1}{z-1}
\end{aligned}
\end{equation}
as $x\rightarrow 0$, so $B_{1}(u,z)=\frac{1}{z-1}$ and $B_{1}(z)=\frac{1}{z-1}.$
\end{remark}

\begin{remark} Eq. (\ref{APN-Stirl-form}) recovers a formula by Apostol (see \cite[Eq. (3.7)]{Ap}), and recently, Xu and Chen~\cite{XC} provided another formula for the Apostol-Bernoulli numbers as follows,  $$B_n(z)=(-1)^{n-1}n\sum_{k=1}^n\frac{(k-1)!}{(z-1)^k}S(n,k)$$ (see \cite[Theorem 4.1]{XC}).
\end{remark}

\begin{theorem}\label{th-main2}
Suppose that $z\not=1$. Under the conventions that $\binom{0}{0}=1$ and $\binom{p}{q}=0$ for $q>p\ge0$, the Apostol-Bernoulli polynomials
$B_{n+1}(u,z)$ for $n\in\mathbb N$ may be expressed as
\begin{equation}\label{ABern-Polyn-determ}
B_{n+1}(u,z)=\frac{(-1)^{n}(n+1)}{(z-1)^{n+1}}\left|\binom{\ell}{m}\left[z(1-u)^{\ell-m}-(-u)^{\ell-m}\right]\right|_{1\le \ell\le n,0\le m\le n-1},
\end{equation}
where $|\cdot|_{1\le \ell\le n,0\le m\le n-1}$ denotes a $n\times n$ determinant.
Consequently, the Apostol-Bernoulli numbers $B_{n+1}(z)$ for $n\in\mathbb N$ can be represented as
\begin{equation}\label{ABern-No-determ}
B_{n+1}(z)=\frac{(-1)^{n}(n+1)}{(z-1)^{n+1}}\left|\binom{\ell}{m}(z-\delta_{\ell m})\right|_{1\le \ell\le n,0\le m\le n-1},
\end{equation}
where the Kronecker delta $\delta_{\ell m}$ is 1 if the variables are equal, and 0 otherwise.
\end{theorem}

\section{Bell polynomials}
As in \cite{Qi}, our proofs are also based on following properties of Bell polynomials.

The Bell polynomials of the second kind $\bell_{n,k}(x_1,x_2,\dotsc,x_{n-k+1})$ are defined by
\begin{multline}\label{Bd}
\bell_{n,k}(x_1,x_2,\dotsc,x_{n-k+1}) \\
=\sum
\frac{n!}{\ell_1!\cdots\ell_{n-k+1}!}\left(\frac{x_1}{1!}\right)^{\ell_1}\left(\frac{x_2}{2!}\right)^{\ell_2}\cdots
\left(\frac{x_{n-k+1}}{{n-k+1}!}\right)^{\ell_{n-k+1}},
\end{multline}
where the sum is taken over all sequences $\ell_1,\ldots,\ell_{n-k+1}$ of non-negative integers such that
$$\ell_1+\ldots+\ell_{n-k+1}=k\quad \text{ and }\quad \ell_1+2\ell_2+\cdots+(n-k+1)\ell_{n-k+1}=n$$
for $n\ge k\ge0$. See~\cite[p.~134, Theorem~A]{Comtet-Combinatorics-74}.

\begin{lemma}[{\cite[Lemma 2.1]{Qi}}]
For $n\ge k\ge0$, the Bell polynomials of the second kind $\bell_{n,k}$ meets
\begin{multline}\label{Bell-times-eq}
\bell_{n,k}(x_1+y_1,x_2+y_2,\dotsc,x_{n-k+1}+y_{n-k+1})\\
=\sum_{r+s=k}\sum_{\ell+m=n}\binom{n}{\ell}\bell_{\ell,r}(x_1,x_2,\dotsc,x_{\ell-r+1})
\bell_{m,s}(y_1,y_2,\dotsc,y_{m-s+1}).
\end{multline}
\end{lemma}

\begin{lemma}[{\cite[p.~135]{Comtet-Combinatorics-74}, ~\cite[Lemma 2.2]{Qi}}]
For $n\ge k\ge0$, we have
\begin{equation}\label{Bell(n-k)}
\bell_{n,k}\left(abx_1,ab^2x_2,\dotsc,ab^{n-k+1}x_{n-k+1}\right) =a^kb^n\bell_{n,k}(x_1,x_2,\dotsc,x_{n-k+1}),
\end{equation}
where $a$ and $b$ are any complex numbers.
\end{lemma}

\begin{lemma}[{\cite[p.~135, Theorem B,~{[3g]}]{Comtet-Combinatorics-74}}]For $n\ge k\ge0$, we have
\begin{equation}\label{Bell3}\bell_{n,k}(\underbrace{1,1,\ldots,1}_{n-k+1})=S(n,k).\end{equation}
\end{lemma}

\section{Proof of Theorem~\ref{th-main1}}
Set
$$m(x)=x,\quad g(x)=\frac{z e^x-1}{e^{ux}}\quad\text{and}\quad f(y)=\frac{1}{y},$$
then we have
$$h(x)=m(x) (f\circ g)(x)=\frac{xe^{ux}}{z e^x-1}$$
which is the generating function of the Apostol-Bernoulli polynomials ~(\ref{AB-def}).
Thus by (\ref{AB-def}), we have
\begin{equation}\label{P1}
\begin{aligned}
\frac{\td^{n+1}}{\td x^{n+1}}h(x)\biggl|_{x=0}
&=B_{n+1}(u,z).
\end{aligned}
\end{equation}
On the other hand, since
\begin{equation}
\begin{aligned}
g(x)&=\frac{z e^x-1}{e^{ux}}=ze^{(1-u)x}-e^{-ux}\\
& = z\sum_{m=0}^{\infty}\frac{(1-u)^{m}}{m!}x^{m}-\sum_{m=0}^{\infty}\frac{(-u)^{m}}{m!}x^{m}\\
&\rightarrow z-1
\end{aligned}
\end{equation}
and
\begin{equation}\label{P1-1}
\begin{aligned}
g'(x)&=z\sum_{m=1}^{\infty}\frac{(1-u)^{m}}{(m-1)!}x^{m-1}-\sum_{m=1}^{\infty}\frac{(-u)^{m}}{(m-1)!} x^{m-1}\\
&\rightarrow z(1-u)-(-u),\\
 g''(x)&=z\sum_{m=2}^{\infty}\frac{(1-u)^{m}}{(m-2)!}x^{m-2}-\sum_{m=2}^{\infty}\frac{(-u)^{m}}{(m-2)!}x^{m-2} \\
&\rightarrow z(1-u)^{2}-(-u)^{2} ,\\
&\;\,\vdots\\
g^{(n-k+1)}(x)&=z\sum_{m=n-k+1}^{\infty}\frac{(1-u)^{m}}{(m-(n-k+1))!}x^{m-(n-k+1)} \\
&\quad-\sum_{m=n-k+1}^{\infty}\frac{(-u)^{m}}{(m-(n-k+1)!}x^{m-(n-k+1)} \\
&\rightarrow z(1-u)^{n-k+1}-(-u)^{n-k+1}
\end{aligned}
\end{equation}
as $x\rightarrow 0.$

In terms of the Bell polynomials of the second kind $B_{n,k},$ the  Fa\`a di Bruno formula for computing higher order derivatives of
composite functions is described in \cite[p.~139, Theorem C]{Comtet-Combinatorics-74} by
\begin{equation}\label{fabr}
\frac{\td^n}{\td x^n}(f\circ g)(x)=\sum_{k=0}^nf^{(k)}(g(x)) \bell_{n,k}\left(g'(x),g''(x),\dotsc,g^{(n-k+1)}(x)\right)
\end{equation}
(see also \cite[p.~93, (3.1)]{Qi}).

By the integral expression (\ref{AB-def}) and (\ref{P1-1}), applying the formula (\ref{fabr}) to the functions
$$f(y)=\frac1y \quad\text{and}\quad y=g(x)=\frac{ze^x-1}{e^{ux}},$$
we have
\begin{equation}\label{foru}
\begin{aligned}
\frac{\td^n}{\td x^n}(f\circ g)(x)&=\sum_{k=0}^nf^{(k)}(g(x)) \bell_{n,k}\left(g'(x),g''(x),\dotsc,g^{(n-k+1)}(x)\right)\\
&=\sum_{k=0}^{n}\frac{(-1)^{k}k!}{(g(x))^{k+1}}\bell_{n,k}\left(g'(x),g''(x),\dotsc,g^{(n-k+1)}(x)\right)\\
&\rightarrow \sum_{k=0}^{n}\frac{(-1)^{k}k!}{(z-1)^{k+1}}\\
&\quad \times\bell_{n,k}\left(z(1-u)-(-u),z(1-u)^{2}-(-u)^{2} ,\dots,z(1-u)^{n-k+1}-(-u)^{n-k+1}\right)
\end{aligned}
\end{equation}
as $x\rightarrow 0$ and
\begin{equation}\label{P2}
\begin{aligned}
\frac{\td^n}{\td x^n}(f\circ g)(x)\biggl|_{x=0}
&=\sum_{k=0}^n\frac{(-1)^{k}k!}{(z-1)^{k+1}} \sum_{r+s=k}\sum_{\ell+m=n}\binom{n}{\ell}\\
&\quad\times\bell_{\ell,r}\left(z(1-u), z(1-u)^{2}, \dotsc,z(1-u)^{\ell-r+1}\right)\\
&\quad\times\bell_{m,s}\left(-(-u),-(-u)^{2}, \dotsc,-(-u)^{m-s+1}\right)\\
&\quad\text{(by \eqref{Bell-times-eq})}
 \\&=\sum_{k=0}^n\frac{(-1)^{k}k!}{(z-1)^{k+1}}\sum_{r+s=k}\sum_{\ell+m=n}\binom{n}{\ell}\\
&\quad\times z^{r}(1-u)^{\ell} \bell_{\ell,r}(\underbrace{1,1,\ldots,1}_{\ell-r+1})\\
&\quad\times (-1)^{s}(-u)^m \bell_{m,s}(\underbrace{1, 1, \ldots,1}_{m-s+1})\\
&\quad\text{(by \eqref{Bell(n-k)})}
\\&=\sum_{k=0}^n\frac{(-1)^{k}k!}{(z-1)^{k+1}}\sum_{r+s=k}\sum_{\ell+m=n}\binom{n}{\ell}\\
&\quad\times z^{r}(1-u)^{\ell}S(\ell,r)
 (-1)^{s}(-u)^m S(m,s)\\
&\quad\text{(by \eqref{Bell3})}.
\end{aligned}
\end{equation}
Thus by Lebnitz's formula for the $n$th derivative of the product of two functions (see~\cite[p.~210, Example 24]{Zorich}), we have
\begin{equation}\label{PR}
\begin{aligned}
\frac{\td^n}{\td x^n}h(x)&=\frac{\td^n}{\td x^n}m(x)(f\circ g)(x)\\
&=\sum_{i=0}^{n}\binom{n}{i}(f\circ g)(x)^{(n-i)}m(x)^{(i)}\\
&=x\frac{\td^n}{\td x^n}(f\circ g)(x)+n\frac{\td^{n-1}}{\td x^{n-1}}(f\circ g)(x)\\
&\quad(\text{since}~m(x)=x)\\
&\rightarrow n\sum_{k=0}^{n-1}\frac{(-1)^{k}k!}{(z-1)^{k+1}}\sum_{r+s=k}\sum_{\ell+m=n-1}\binom{n-1}{\ell}\\
&\quad\times z^{r}(1-u)^{\ell}S(\ell,r)(-1)^s(-u)^m S(m,s)\\
&\quad\text{(by \eqref{P2})}
\end{aligned}
\end{equation}
as $x\rightarrow 0$.
Then by comparing (\ref{P1}) with (\ref{PR}), we get our result.

Finally, letting $u=0$ in (\ref{AP-Stirl-form}), we obtain (\ref{APN-Stirl-form}).

\section{Proof of Theorem~\ref{th-main2}}

Let $\mu=\mu(x)$ and $\nu=\nu(x)\neq0$ be differentiable functions.
Set
$$\frac{\td^n}{\td x^n}\left(\frac{\mu}{\nu}\right)=\frac{(-1)^nw_n}{\nu^{n+1}}$$
at every point $\nu(x)\neq0.$
By \cite[p.~40]{Bour}, we have
\begin{equation}\label{matrix}
w_n=
 \left|
   \begin{array}{cccccc}
     \mu         & \nu         & 0                       & 0                       & \cdots & 0                 \\
     \mu'        & \nu'        & \nu                     & 0                       & \cdots & 0                 \\
     \mu''       & \nu''       & 2\nu'                   & \nu                     & \cdots & 0                 \\
     \vdots      & \vdots      & \vdots                  & \vdots                  & \vdots & \vdots            \\
     \mu^{(n-1)} & \nu^{(n-1)} & \binom{n-1}1\nu^{(n-2)} & \binom{n-1}2\nu^{(n-3)} & \cdots & \nu               \\
     \mu^{(n)}   & \nu^{(n)}   & \binom{n}1\nu^{(n-1)}   & \binom{n}2\nu^{(n-2)}   & \cdots & \binom n{n-1}\nu' \\
   \end{array}
 \right|.
\end{equation}
As in \cite[p. 94--95, the~ first~ proof~ of~ Theorem 1.2]{Qi},
we may reformulate the formula (\ref{matrix}) as
\begin{equation}\label{reform-matrix}
\frac{\td^n}{\td x^n}\left(\frac{\mu}{\nu}\right)
=\frac{(-1)^n}{\nu^{n+1}}\begin{vmatrix} A_{(n+1)\times 1} & B_{(n+1)\times n}\end{vmatrix}_{(n+1)\times(n+1)} ,
\end{equation}
where the matrices
$$A_{(n+1)\times 1}=(a_{\ell,1})_{0\leq\ell\leq n}=
\left(
   \begin{array}{c}
     \mu         \\
     \mu'        \\
     \mu''       \\
     \vdots      \\
     \mu^{(n-1)} \\
     \mu^{(n)}   \\
   \end{array}
 \right)$$
and
$$\begin{aligned}
B_{(n+1)\times n}&=(b_{\ell,m})_{0\leq \ell\leq n,0\leq m\leq n-1} \\
&=
\left(
   \begin{array}{cccccc}
     \nu         & 0                       & 0                       & \cdots & 0                 \\
     \nu'        & \nu                     & 0                       & \cdots & 0                 \\
     \nu''       & 2\nu'                   & \nu                     & \cdots & 0                 \\
     \vdots      & \vdots                  & \vdots                  & \vdots & \vdots            \\
     \nu^{(n-1)} & \binom{n-1}1\nu^{(n-2)} & \binom{n-1}2\nu^{(n-3)} & \cdots & \nu               \\
     \nu^{(n)}   & \binom{n}1\nu^{(n-1)}   & \binom{n}2\nu^{(n-2)}   & \cdots & \binom n{n-1}\nu' \\
   \end{array}
 \right)
\end{aligned}
$$
satisfy
$$a_{\ell,1}=\mu^{(\ell)}(x)\quad\text{and}\quad b_{\ell,m}=\binom \ell m \nu^{(\ell-m)}(x)$$
under the conventions that $\nu^{(0)}(x)=\nu(x)$ and that $\binom \ell m=0$ and $\nu^{(\ell-m)}(x)\equiv0$ for $\ell< m$
(see \cite[(3.2) and (3.3)]{Qi}).
If we let
\begin{equation}
\mu(x)=1, \quad  \nu(x)=\frac{ze^{x}-1}{e^{ux}},
\end{equation}
then in this case,
the $a_{\ell,1}$ becomes
\begin{equation}\label{a ell}
a_{0,1}=1,\quad a_{\ell,1}=0 \quad\text{for }\ell\geq1
\end{equation}
and the $b_{\ell,m}$ becomes
\begin{equation}\label{P3}
\begin{aligned}
\binom{\ell}{m}\nu^{(\ell-m)}(x)&=\binom{\ell}{m}\frac{\td^{\ell-m}}{\td x^{\ell-m}}\left(\frac{ze^{x}-1}{e^{ux}}\right)\\
&=\binom{\ell}{m}\frac{\td^{\ell-m}}{\td x^{\ell-m}}\left(ze^{(1-u)x}-e^{-ux}\right)\\
&=\binom{\ell}{m}\frac{\td^{\ell-m}}{\td x^{\ell-m}}\left(z\sum_{k=0}^{\infty}\frac{(1-u)^{k}x^{k}}{k!}-\sum_{k=0}^{\infty}\frac{(-u)^{k}x^{k}}{k!}\right)\\
&=\binom{\ell}{m}\left(z\sum_{k=\ell-m}^{\infty}\frac{(1-u)^{k}x^{k-(\ell-m)}}{(k-(\ell-m))!}
-\sum_{k=l-m}^{\infty}\frac{(-u)^{k}x^{k-(\ell-m)}}{(k-(\ell-m))!}\right)\\
&\rightarrow\binom{\ell}{m} \left[z(1-u)^{\ell-m}-(-u)^{\ell-m}\right],\quad x\rightarrow 0
\end{aligned}
\end{equation}
for $0\leq\ell\leq n$ and $0\leq m\leq n-1$ with $\ell\geq m.$
Thus,
as in \cite[p. 95, the~ first~ proof~ of~ Theorem 1.2]{Qi}, by (\ref{reform-matrix}), (\ref{a ell}) and (\ref{P3}), we have
\begin{equation}\label{PRth2}
\begin{aligned}
\frac{\td^n}{\td x^n}\left(\frac{\mu}{\nu}\right)
&=\frac{(-1)^n}{\nu^{n+1}}
\left|
   \begin{array}{cccccc}
     1         & b_{0,0}         & 0                     & 0                     & \cdots & 0 \\
     0         & b_{1,0}         & b_{1,1}               & 0                     & \cdots & 0 \\
     0         & b_{2,0}         & b_{2,1}               & b_{2,2}               & \cdots & 0 \\
     \vdots    & \vdots          & \vdots                & \vdots                & \vdots & \vdots \\
     0         & b_{n-1,0}       & b_{n-1,1}             & b_{n-1,2}             & \cdots & b_{n-1,n-1} \\
     0         & b_{n,0}         & b_{n,1}               & b_{n,2}               & \cdots & b_{n,n-1} \\
   \end{array}
 \right| \\
&=\frac{(-1)^n}{\nu^{n+1}} \left|b_{\ell,m}\right|_{1\le \ell\le n,0\le m\le n-1}\\
&\to \frac{(-1)^n}{(z-1)^{n+1}}\left|\binom{\ell}{m}\left[z(1-u)^{\ell-m}-(-u)^{\ell-m}\right]\right|_{1\le \ell\le n,0\le m\le n-1}
\end{aligned}
\end{equation}
as $x\rightarrow 0.$
Then again by Lebnitz's formula for the $n$-th derivative of the product of two functions (see~\cite[p.~210, Example 24]{Zorich}), we have
\begin{equation}\label{PR2}
\begin{aligned}
\frac{\td^{n+1}}{\td x^{n+1}}\left(\frac{xe^{ux}}{ze^{x}-1}\right)&=\frac{\td^{n+1}}{\td x^{n+1}}\left(x\left(\frac{\mu}{\nu}\right)\right)\\
&=\sum_{i=0}^{n+1}\binom{n+1}{i}\left(\frac{\mu}{\nu}\right)^{(n+1-i)}x^{(i)}\\
&=x\frac{\td^{n+1}}{\td x^{n+1}}\left(\frac{\mu}{\nu}\right)+(n+1)\frac{\td^{n}}{\td x^{n}}\left(\frac{\mu}{\nu}\right)\\
&\rightarrow (n+1)\frac{\td^{n}}{\td x^{n}}\left(\frac{\mu}{\nu}\right) \\
&=\frac{(-1)^{n}(n+1)}{(z-1)^{n+1}}\left|\binom{\ell}{m}\left[z(1-u)^{\ell-m}-(-u)^{\ell-m}\right]\right|_{1\le \ell\le n,0\le m\le n-1} \\
&\quad\text{(by \eqref{PRth2})}
\end{aligned}
\end{equation}
as $x\rightarrow 0$.
Thus, by comparing (\ref{P1}) with (\ref{PR2}), we obtain (\ref{ABern-Polyn-determ}).

Finally, letting $u=0$ in (\ref{ABern-Polyn-determ}), we obtain (\ref{ABern-No-determ}).

\section*{Acknowledgment}

The authors are enormously grateful to the anonymous referee whose comments and suggestions lead to a large improvement of the
paper.

\bibliography{central}

\end{document}